\documentclass[a4paper,12pt]{article}
\usepackage{amsmath}
\usepackage{amsfonts}
\usepackage{amssymb}
\usepackage{latexsym}
\usepackage{epsfig}
\usepackage{graphicx}
\usepackage{oldgerm}
\usepackage{theorem}

\def\R{\mathbb{R}}
\def\C{\mathbb{C}}

\def\W{\mathbb{W}}
\def\GT{\mathbb{GT}}
\def\bE{{\bf E}}

\def\1{{\bf 1}}

\def\cN{{\cal N}}
\def\cL{{\cal L}}
\def\cH{{\cal H}}

\def\cS{{\cal S}}

\def\x{\mib{x}}
\def\y{\mib{y}}
\def\z{\mib{z}}
\def\B{\mib{B}}
\def\k{\mib{k}}
\def\Z{\mib{Z}}
\def\T{\mib{T}}
\def\t{\mib{t}}

\def\vnu{\mib{\nu}}
\def\vmu{\mib{\mu}}
\def\veta{\mib{\eta}}
\def\valpha{\mib{\alpha}}

\def\gl{\mathfrak{gl}}

\theorembodyfont{\itshape}

\newtheorem{thm}{Theorem}[section]
\newtheorem{lem}[thm]{Lemma}

\newtheorem{prop}[thm]{Proposition}

\newcommand{\mib}[1]{\mbox{\boldmath $#1$}}
\newcommand{\SSC}[1]{\section{#1}\setcounter{equation}{0}}
\newcommand{\qed}{\hbox{\rule[-2pt]{3pt}{6pt}}}

\begin{document}

\title{\bf Survival probability of mutually killing Brownian motions
and the O'Connell process}
\author{
Makoto Katori
\footnote{
Department of Physics,
Faculty of Science and Engineering,
Chuo University, 
Kasuga, Bunkyo-ku, Tokyo 112-8551, Japan;
e-mail: katori@phys.chuo-u.ac.jp
}}
\date{23 March 2012}
\pagestyle{plain}
\maketitle
\begin{abstract}
Recently O'Connell introduced an interacting diffusive particle
system in order to study a directed polymer model in 1+1 dimensions.
The infinitesimal generator of the process 
is a harmonic transform of the quantum Toda-lattice Hamiltonian
by the Whittaker function.
As a physical interpretation of this construction,
we show that the O'Connell process without drift
is realized as a system of mutually killing Brownian motions
conditioned that all particles survive forever.
When the characteristic length of interaction
killing other particles goes to zero, the process is reduced to
the noncolliding Brownian motion (the Dyson model).

\end{abstract}


\SSC{Introduction}

In this paper we introduce a system of finite number of
one-dimensional Brownian motions which kill each other,
and evaluate long-term asymptotics of
the probability that all particles survive. 
Then we define a process of 
mutually killing Brownian motions 
{\it conditioned that all particles survive forever}.
We show that this conditional process is equivalent to
a special case of the process recently introduced by
O'Connell in order to analyze a directed polymer model
in 1+1 dimensions \cite{OCo09}.

\begin{figure}
\includegraphics[width=0.8\linewidth]{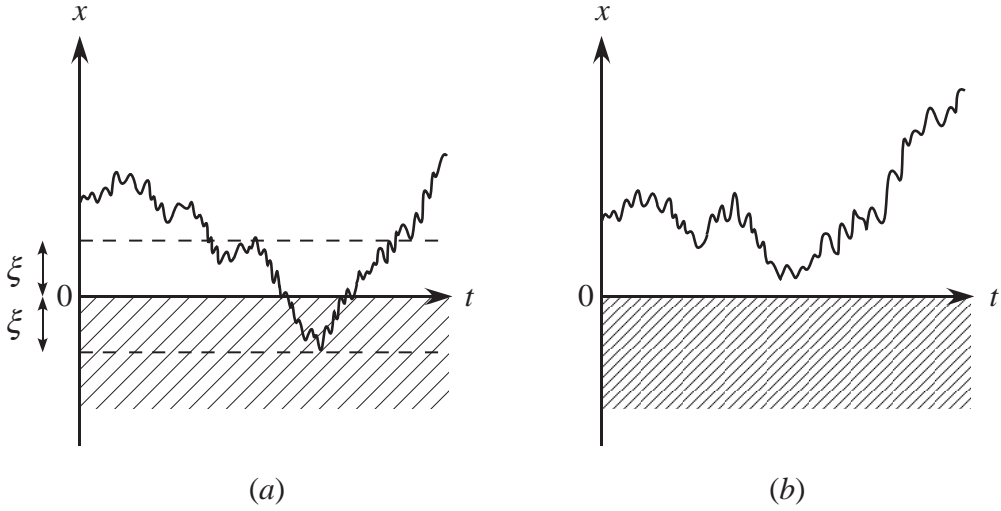}
\caption{
(a) An illustration of a path of surviving Brownian particle
in one-dimension with a killing term
$-V(x)$ given by (\ref{eqn:V1}).
The parameter $\xi$ is the characteristic length of the interaction
to kill a particle acting from the origin.
It also indicates the penetration length of a particle
into the negative region.
(b) An illustration of a path of surviving Brownian 
particle in the absorbing Brownian motion 
with an absorbing wall at the origin.
It is regarded as the limit $\xi \to 0$ of the
system shown in (a).
\label{fig:Fig1}}
\end{figure}

As an introduction of our study of many-particle systems,
here we consider a family of one-particle systems
with a parameter $\xi > 0$.
Let $B(t), t \geq 0$ be such a one-dimensional Brownian motion
that its survival probability $P(t)$ decays following the
equation
\begin{equation}
\frac{d P(t)}{dt}=-V(B(t)) P(t),
\quad t \geq 0,
\label{eqn:Pt1}
\end{equation}
with a decay rate function
\begin{equation}
V(x)=\frac{1}{2 \xi^2} e^{-2x/\xi}, \quad x \in \R.
\label{eqn:V1}
\end{equation}
The function (\ref{eqn:V1}) implies that, 
if the Brownian particle moves in the positive region
far from the origin $x \gg \xi$, 
decay of survival probability is 
negligible, while as it approaches the
origin the decay becomes large.
Note that the Brownian particle is able to 
penetrate a negative region $x < 0$,
but there the decay rate of survival probability
grows exponentially as a function of $|x|$.
The parameter $\xi$ is the {\it characteristic length
of the interaction} to kill a particle
acting from the origin and it is also the 
{\it penetration length} of a particle 
into the negative region.
See Fig.\ref{fig:Fig1}(a).
For $0 < t < \infty$, if a path of the Brownian motion
up to time $t$ is given as $\{B(s): 0 \leq s \leq t\}$,
the survival probability of the particle at the time $t$
is given by an integration of (\ref{eqn:Pt1}), 
\begin{eqnarray}
P(t | \{B(s): 0 \leq s \leq t\})
&=& \exp \left\{-\int_0^t V(B(s)) ds \right\}
\nonumber\\
&=& \exp \left\{
-\frac{1}{2 \xi^2} \int_0^{t} 
e^{-2B(s)/\xi} ds \right\}.
\label{eqn:Pt2}
\end{eqnarray}
Provided that $B(0)>0$, in the limit
$\xi \to 0$, the process becomes
an absorbing Brownian motion in the positive region 
with an absorbing wall
at the origin, in which the survival probability
(\ref{eqn:Pt2}) is zero if the particle hits
the origin at any time $s \leq t$, and it is one
if it stays in the positive region 
in the time period $[0, t]$. See Fig.\ref{fig:Fig1}(b).

Let $Q(t, y|x)$ be the transition probability density
of this killing Brownian motion from a position $x \in \R$
to a position $y \in \R$ during time $t >0$.
It will be obtained by taking an average of (\ref{eqn:Pt2})
over all realizations of Brownian path under the conditions
$B(0)=y, B(t)=x$. 
(Note that the time ordering is reversed,
since the equation of $Q$ given by (\ref{eqn:Q2}) below
is the {\it backward} Kolmogorov equation.)
That is,
\begin{equation}
Q(t, y|x)= \bE \left[
\exp \left\{ -\int_0^t V(B(s)) ds \right\} 
\1(B(0)=y, B(t)=x) \right],
\label{eqn:Q1}
\end{equation}
where $\bE[\ \cdot \ ]$ denotes an expectation
with respect to Brownian motion, and
$\1(\omega)$ is the indicator function of
a condition $\omega$; $\1(\omega)=1$ if
$\omega$ is satisfied, and $\1(\omega)=0$ otherwise.
It is the Feynman-Kac formula (see, {\it e.g.}, \cite{KS91})
for the unique solution of the diffusion equation
\begin{equation}
\frac{\partial}{\partial t} Q(t, y|x)
=\left( \frac{1}{2} \frac{\partial^2}{\partial x^2}
-V(x) \right) Q(t, y|x)
\label{eqn:Q2}
\end{equation}
under the initial condition
$Q(0, y|x)=\delta(x-y)$.

Let $I_{\nu}(z), z, \nu \in \C$ be the modified Bessel function
of the first kind
\begin{equation}
I_{\nu}(z)=\sum_{k=0}^{\infty} 
\frac{(z/2)^{\nu+2k}}{\Gamma(k+1) \Gamma(k+\nu+1)},
\quad |z| < \infty, \quad
|{\rm arg} \, z| < \pi,
\label{eqn:Inu1}
\end{equation}
and $K_{\nu}(z), z \in \C$ be Macdonald's function \cite{Wat44,Leb65}
defined by
\begin{equation}
K_{\nu}(z)=\frac{\pi}{2}
\frac{I_{-\nu}(z)-I_{\nu}(z)}{\sin (\nu \pi)},
\quad |{\rm arg} \, z|< \pi, \quad \mbox{for} \quad
\nu \not= 0, \pm 1, \pm 2, \dots,
\label{eqn:Knu1}
\end{equation}
and for integers $\nu=n$,
\begin{equation}
K_{n}(z)= \lim_{\nu \to n} K_{\nu}(z),
\quad n=0, \pm 1, \pm 2, \dots,
\label{eqn:Knu2}
\end{equation}
which are both analytic functions
of $z$ for all $z$ in the complex
plane $\C$ cut along the negative real axis,
and entire functions of $\nu$.
We see that $I_{\nu}(z)$ and $K_{\nu}(z)$
are linearly independent solutions
of the differential equation
$$
\frac{d^2 w}{dz^2}
+\frac{1}{z} \frac{dw}{dz}
-\left( 1+ \frac{\nu^2}{z^2} \right) w=0.
$$
For $x >0$ and $\nu \geq 0$,
$I_{\nu}(x)$ is a positive function
which increases monotonically as $x \to \infty$,
while $K_{\nu}(x)$ is a positive function
which decreases monotonically
as $x \to \infty$.
Then an integral representation of $Q(t,y|x)$ is given
by \cite{MY00,MY05,OCo09,Kat11}
\begin{equation}
Q(t, y|x)
= \frac{1}{\pi} \int_{-\infty}^{\infty}
e^{-k^2 t/2} K_{i \xi k}(e^{-x/\xi})
K_{-i \xi k}(e^{-y/\xi}) 
|\Gamma(i \xi k)|^{-2} d k,
\label{eqn:Q3}
\end{equation}
where $i=\sqrt{-1}$ and $\Gamma(z)$ is the gamma function.

For each initial position $x \in \R$, the survival probability
of this killing Brownian motion at time $T \geq 0$
is obtained by integrating $Q(T,y|x)$ over $y \in \R$,
\begin{equation}
\cN(T,x)=\int_{-\infty}^{\infty} Q(T,y|x) dy,
\quad x \in \R, \quad T \geq 0.
\label{eqn:N1}
\end{equation}
If we condition the process to survive up to time $T>0$, 
the transition probability density
from the state $(s,x)$ to $(t,y)$ of the
killing Brownian motion is given by
\begin{equation}
P^{T}(s,x;t,y)
=\frac{\cN(T-t,y)}{\cN(T-s,x)} Q(t-s, y|x),
\quad x, y \in \R, \quad 0 \leq s \leq t \leq T.
\label{eqn:PT1}
\end{equation}
If we use the integral representations of
$K_{\nu}(z)$ \cite{Wat44,Leb65},
\begin{equation}
K_{\nu}(z)=\frac{2^{\nu} \Gamma(\nu+1/2)}{z^{\nu} \sqrt{\pi}}
\int_0^{\infty} \frac{\cos (z u)}{(1+u^2)^{\nu+1/2}} du,
\label{eqn:KInt1}
\end{equation}
and
\begin{equation}
K_{\nu}(z)=\frac{1}{2} \int_0^{\infty} s^{\nu-1}
\exp \left\{-\frac{z}{2} \left(s+\frac{1}{s} \right) \right\} ds,
\label{eqn:KInt2}
\end{equation}
we can evaluate the long-term asymptotics of the survival probability
for $|x| < \infty$, 
\begin{equation}
\cN(T,x) = C T^{-1/2} K_0(e^{-x/\xi})
\times \left\{ 1+{\cal O} \left( \frac{x}{\sqrt{T}} \right) \right\}
\quad \mbox{in} \quad \frac{x}{\sqrt{T}} \to 0
\label{eqn:Nlimit}
\end{equation}
with $C=3 \xi \sqrt{2/\pi}$ as shown in Appendix A.
Then we can take the limit $T \to \infty$ of (\ref{eqn:PT1}) as
\begin{eqnarray}
P(t,y|x) &\equiv& \lim_{T \to \infty} P^{T}(0,x;t,y)
\nonumber\\
&=& \frac{K_0(e^{-y/\xi})}{K_0(e^{-x/\xi})}
Q(t, y|x), \quad x, y \in \R, \quad 0 \leq t < \infty.
\label{eqn:P1}
\end{eqnarray}
This is the transition probability density of the
killing Brownian motion {\it conditioned to survive forever}.
We should note that this conditional Brownian motion
is the diffusion process with infinitesimal generator
\begin{equation}
\cL_{\rm MY}=\frac{1}{2} \frac{d^2}{dx^2}
+\frac{d}{dx} \Big\{ \log K_0(e^{-x/\xi}) \Big\} \frac{d}{dx}.
\label{eqn:MY1}
\end{equation}
Matsumoto and Yor \cite{MY00,MY05} have studied the exponential
functionals of Brownian motion
\begin{equation}
Z_{\rm MY}(t)= \xi \log \left\{ \frac{1}{\xi^2}
\int_0^t e^{2 B(s)/\xi} ds \right\} - B(t), \quad
t \geq 0,
\label{eqn:MY2}
\end{equation}
which is the diffusion process whose infinitesimal
generator is (\ref{eqn:MY1}).
We can say that the Matsumoto-Yor process (\ref{eqn:MY2})
is realized as the present killing Brownian motion
conditioned to survive forever \cite{Kat11}.

Now we consider the limit $\xi \to 0$ of the formulas
obtained above. Since $I_{\nu}(z)$ is defined as
the series expansion (\ref{eqn:Inu1}), we can see
$I_{i \xi k}(e^{-x/\xi}) \simeq (e^{-x/\xi}/2)^{i \xi k}/\Gamma(i\xi k+1)
\simeq e^{-i k x}$ in $\xi \to 0$ for $x > 0$.
Then (\ref{eqn:Knu1}) gives
$K_{i \xi k}(e^{-x/\xi}) \simeq \sin(kx)/(\xi k) \1(x >0)$
in $\xi \to 0$.
Therefore the integral of (\ref{eqn:Q3}) can be performed
in the limit and we have
\begin{eqnarray}
&& q(t,y|x) \equiv \lim_{\xi \to 0} Q(t,y|x)
\nonumber\\
&& \qquad = \frac{1}{\sqrt{2 \pi t}}
\Big\{ e^{-(x-y)^2/2t} - e^{-(x+y)^2/2t} \Big\}
\, \1(x>0, y \geq 0).
\label{eqn:q1}
\end{eqnarray}
It is the transition probability density of the
absorbing Brownian motion with an absorbing wall at the origin,
which is easily obtained by applying the
reflection principle of Brownian motion \cite{KS91}.
Applying the result (\ref{eqn:K0limA}) 
shown in Appendix A, we see that
\begin{equation}
\lim_{\xi \to 0} \xi K_0(e^{-x/\xi}) 
=x \, \1(x>0), 
\label{eqn:KlimB}
\end{equation}
and then (\ref{eqn:P1}) gives in this limit
\begin{eqnarray}
p(t,y|x) &\equiv& \lim_{\xi \to 0} P(t,y|x)
\nonumber\\
&=& \frac{y}{x} q(t,y|x), \quad
x>0, y \geq 0, \quad t \geq 0,
\label{eqn:p1}
\end{eqnarray}
which is a harmonic transform ($h$-transform) of (\ref{eqn:q1})
and is identified with the transition probability density
of the three-dimensional Bessel process, BES(3)
(see, for instance, \cite{KT_Sugaku_11}).
As a matter of fact \cite{MY00,MY05},
(\ref{eqn:MY2}) gives
$\lim_{\xi \to 0} Z_{\rm MY}(t) =
2 \sup_{0 \leq s \leq t} B(s)-B(t), t \geq 0$,
which is indeed distributed as BES(3)
(Pitman's $2M-X$ theorem \cite{Pit75}).

In the present paper, we consider an $N$-particle system
of one-dimensional Brownian motions with $N \geq 2$,
$\B(t)=(B_1(t), \dots, B_N(t)), t \geq 0$, 
with a positive parameter $\xi >0$, 
such that
the probability that all $N$ particles
survive up to time $t$, $P_N(t)$, 
decays following the equation
\begin{equation}
\frac{dP_N(t)}{dt}=-V_N(\B(t)) P_N(t),
\quad t \geq 0
\label{eqn:PN1}
\end{equation}
with a decay rate function
\begin{equation}
V_N(\x)=\frac{1}{\xi^2} \sum_{j=1}^{N-1} e^{-(x_{j+1}-x_j)/\xi},
\quad \x \in \R^N.
\label{eqn:VN1}
\end{equation}
We study an integral representation of the
transition probability density $Q_N(t,\y|\x), 
\x, \y \in \R^N, t \geq 0$,
which is a multivariate extension of (\ref{eqn:Q3}).
The extensions of the formulas
(\ref{eqn:Nlimit}) and (\ref{eqn:P1}) are shown.
We prove that the $N$-particle system of the mutually killing Brownian motions
following (\ref{eqn:PN1}) and (\ref{eqn:VN1}) 
{\it conditioned that all $N$ particles survive forever} 
is equivalent to
a special case (without drift, $\vnu=0$) of 
the O'Connell process \cite{OCo09,Kat11}.

In the context of quantum mechanics,
the decay rate functions of survival probability
(\ref{eqn:V1}) and (\ref{eqn:VN1}) in the killing 
Brownian motions are considered to giving potential
energy of the systems. They are called 
the {\it Yukawa potential} and the 
{\it Toda-lattice potential} \cite{Toda89}, respectively.
The multivariate extensions of Macdonald's function
(\ref{eqn:Knu1}) with (\ref{eqn:Knu2}) are
the Whittaker functions
(see \cite{BO11,OCo09,COSZ11,BC11,OCo12} and references therein),
which have been extensively studied in mathematical physics as
eigenfunctions of the quantum Toda lattice
\cite{Kos77,Skl85,STS94,Giv97,KL99,KL00,KL01,JK03,GKLO06,GLO08}.
See also \cite{GNSS11}.

In Sect.2, as preliminaries, useful integral representations
of the Whittaker functions are given for the
$GL(N,\R)$-quantum Toda lattice.
The transition probability density of the system of mutually
killing Brownian motions is given,
in the case that killing of particles does not occur
during an observing time period, 
as an integral of
a product of the Whittaker functions over the Sklyanin measure.
Then the main results are given.
We also demonstrate that when the characteristic length $\xi$
of the interaction killing other particles goes to zero,
the O'Connell process is reduced to the noncolliding Brownian motion,
which is known to be equivalent to
the Dyson model \cite{Gra99,KT_Sugaku_11} for the eigenvalue process
of Hermitian matrix-valued diffusion process
({\it i.e.} Dyson's Brownian motion
model with the parameter $\beta=2$ \cite{Dys62}).
The proofs of Lemmas and Proposition are given in Sect.3.
Appendix is given for proving 
(\ref{eqn:Nlimit}) and (\ref{eqn:KlimB}) used above.

\SSC{Preliminaries and Main Results}
\subsection{Quantum Toda lattice and Whittaker functions}

In this subsection  we set $\xi=1$ in Eq.(\ref{eqn:VN1}) 
and consider $V_N(\x)$ as a potential energy of 
a quantum $N$-particle system in one dimension.
Then we have
the Hamiltonian of the $GL(N,\R)$-quantum Toda lattice
\begin{equation}
\cH_N = -\frac{1}{2} \sum_{j=1}^{N} \frac{\partial^2}{\partial x_j^2}
+\sum_{j=1}^{N-1} e^{-(x_{j+1}-x_j)}
\label{eqn:H1}
\end{equation}
for $N \in \{2,3, \dots\}$.
(In this paper, we set the mass of particle $m=1$
and $\hbar =1$ for simplicity of notation.)
The generalized eigenvalue problem
\begin{equation}
\cH_N \Psi_{\lambda}(\x)=\lambda \Psi_{\lambda}(\x),
\quad \x \in \R^N, 
\label{eqn:H2}
\end{equation}
is solved by setting
\begin{equation}
\lambda= - \frac{1}{2} \sum_{j=1}^{N} \nu_j^2
=-\frac{1}{2} |\vnu|^2, \quad
\vnu=(\nu_1, \dots, \nu_N) \in \C^N
\label{eqn:lambda1}
\end{equation}
with the eigenfunction
$\Psi_{\lambda}(\x)=\psi_{\vnu}^{(N)}(\x)$,
where $\psi_{\vnu}^{(N)}(\x)$ is the
$GL(N,\R)$-Whittaker function \cite{Kos77,Skl85,STS94,KL99}
(the class-one $\gl_N$-Whittaker function
\cite{GLO08,BO11,OCo09,COSZ11,BC11}).
When $N=2$, it is expressed by using Macdonald's 
function (\ref{eqn:Knu1}) as
$$
\psi_{(\nu_1, \nu_2)}^{(2)}(x_1,x_2)
=2 e^{(\nu_1+\nu_2)(x_1+x_2)/2}
K_{\nu_2-\nu_1}(2 e^{-(x_2-x_1)/2}).
$$
In this sense, the Whittaker functions 
$\{\psi^{(N)}_{\vnu}(\x)\}_{N \geq 2}$ are regarded as
multivariate extensions of Macdonald's function.

Corresponding to the two kinds of integral
representations (\ref{eqn:KInt1}) and (\ref{eqn:KInt2})
for $K_{\nu}(z)$, two integral representations
are known for the Whittaker function
$\psi_{\vnu}^{(N)}(\x), \vnu \in \C^N$.
The integral expression corresponding to (\ref{eqn:KInt1})
is the classical one originally given by
Jacquet \cite{Jac67} and was rewritten as the
following form by Stade \cite{Sta90}.
(Here we use the notation given by \cite{KL99}.)
Let $\Z$ denote an upper triangular $N \times N$ matrix
with unit diagonal; $\Z=(Z_{j,k})_{1 \leq j, k \leq N}$,
\begin{equation}
Z_{j,k}=\left\{ \begin{array}{ll}
z_{j,k}, & \quad 1 \leq j < k \leq N, \cr
1, & \quad 1 \leq j=k \leq N, \cr
0, & \quad 1 \leq k < j \leq N.
\end{array} \right.
\label{eqn:Z1}
\end{equation}
We write the integral of a function $f$ of $\Z$ over all real $\Z$ as
$$
\int_{\R^{N(N-1)/2}} f(\Z) d \Z
\equiv \prod_{k=1}^{N} \prod_{j=1}^{k-1}
\int_{-\infty}^{\infty} dz_{j,k} f(\Z).
$$
The transpose of $\Z$ is denoted by $^{t}\!\Z$ and the
principal minor of size $j$ of matrix
$\Z \, {^{t}\!\Z}$ is written as
$\Delta_{j}(\Z \, {^{t}\!\Z}), 1 \leq j \leq N$.
For $\x=(x_1, \dots, x_N) \in \R^N, 
\vnu=(\nu_1, \dots, \nu_N) \in \C^N$,
\begin{eqnarray}
&& \psi_{\vnu}^{(N)}(\x)
= \exp(\vnu \cdot \x) \prod_{1 \leq j < k \leq N}
\frac{\Gamma(\nu_k-\nu_j+1/2)}{\sqrt{\pi}}
\nonumber\\
&& \qquad \times \int_{\R^{N(N-1)/2}}
\prod_{j=1}^{N-1}
\Bigg[ \Delta_j(\Z \,{^{t}\!\Z})^{\nu_j-\nu_{j+1}-1/2}
\exp \Big\{2 i z_{j,j+1} e^{-(x_{j+1}-x_j)/2} \Big\} \Bigg] d \Z.
\label{eqn:WInt1}
\end{eqnarray}
Another integral representation, which corresponds to
(\ref{eqn:KInt2}), was given by Givental \cite{Giv97}.
See also \cite{KL01,GKLO06}.
Let $\T$ denote a lower triangular array
with size $N$,
$\T=(T_{j, k}, 1 \leq k \leq j \leq N)$.
For a given $\x \in \R^N$, let $\Gamma_N(\x)$ be
the space of all real triangular arrays $\T$
with size $N$ conditioned
\begin{equation}
T_{N,k}=x_k, \quad 1 \leq k \leq N.
\label{eqn:qToda5}
\end{equation}
We write the integral of a function $f$ of $\T$ 
over $\Gamma_N(\x)$ as
$$
\int_{\Gamma_N(\x)} f(\T) d \T
\equiv \prod_{j=1}^{N} \prod_{k=1}^j
\int_{-\infty}^{\infty} d T_{j, k} \,
f(\T) \prod_{\ell=1}^{N} \delta(T_{N,\ell}-x_{\ell}).
$$
Then 
\begin{eqnarray}
\psi_{\vnu}^{(N)}(\x)
&=&\int_{\Gamma_N(\x)} 
\exp \left[ 
\sum_{j=1}^{N} \nu_{j}
\left( \sum_{k=1}^{j} T_{j, k}
-\sum_{k=1}^{j-1} T_{j-1, k} \right)
\right. \nonumber\\
&& \quad \left.
- \sum_{j=1}^{N-1} \sum_{k=1}^j
\Big\{ e^{-(T_{j,k}-T_{j+1,k})}
+e^{-(T_{j+1, k+1}-T_{j,k})} \Big\} \right] d \T.
\label{eqn:WInt2}
\end{eqnarray}

\vskip 0.3cm
\noindent{\bf Remark 1.} \,
There is the `third' kind of integral representation
of the Whittaker function, 
\begin{eqnarray}
\psi_{\vnu}^{(N)}(\x)
&=& \int_{{\cal S}} \prod_{j=1}^{N-1}
\frac{\prod_{k=1}^j \prod_{\ell=1}^{j+1}
\Gamma(\gamma_{j+1,\ell}-\gamma_{j,k})}
{(2\pi)^j j! \prod_{1 \leq k < \ell \leq j}
|\Gamma(\gamma_{j,\ell}-\gamma_{j,k})|^2}
\nonumber\\
&& \quad \times
\exp \left\{ i \sum_{j=1}^{N} \sum_{k=1}^j
(\gamma_{j,k}-\gamma_{j-1,k}) x_j \right\}
\prod_{j=1}^{N-1} \prod_{k=1}^{j} d \gamma_{j,k},
\label{eqn:WInt3}
\end{eqnarray}
where $(\gamma_{N,1}, \dots, \gamma_{N,N}) 
\equiv (\nu_1, \dots, \nu_N)$,
$\gamma_{j,k}=0$ for $k > j$,
and the domain of integration ${\cal S}$ is
defined by the conditions
$\min_k \{{\rm Im} \gamma_{j,k} \}
> \max_{k}\{{\rm Im} \gamma_{j+1, k} \}$
for all $j=1,2, \dots, N-1$.
This is called the Mellin-Barnes integral representation,
since it can be regarded as the multivariate extension
of the Mellin-Barnes representation of
Macdonald's function \cite{Wat44}
$$
K_{\nu}(z)=\frac{1}{8 \pi i} \int_{{\cal C}}
\Gamma \left( \frac{s+\nu}{2} \right)
\Gamma \left( \frac{s-\nu}{2} \right)
\left( \frac{z}{2} \right)^{-2} ds,
$$
where the path of integration ${\cal C}$ being a vertical line
to the right of any poles of the integrand.
The derivation of the Mellin-Barnes integral representation
(\ref{eqn:WInt3}) from the classical one (\ref{eqn:WInt1})
is shown in \cite{Sta01,IS07} (see also \cite{KL00,KL01}).
The equivalence between the Givental integral representation
(\ref{eqn:WInt2}) and the Mellin-Barnes integral
representation (\ref{eqn:WInt3}) is fully discussed
in \cite{GLO08}.
\vskip 0.3cm

\subsection{Transition probability density of the
mutually killing Brownian motions}

\begin{figure}
\includegraphics[width=0.7\linewidth]{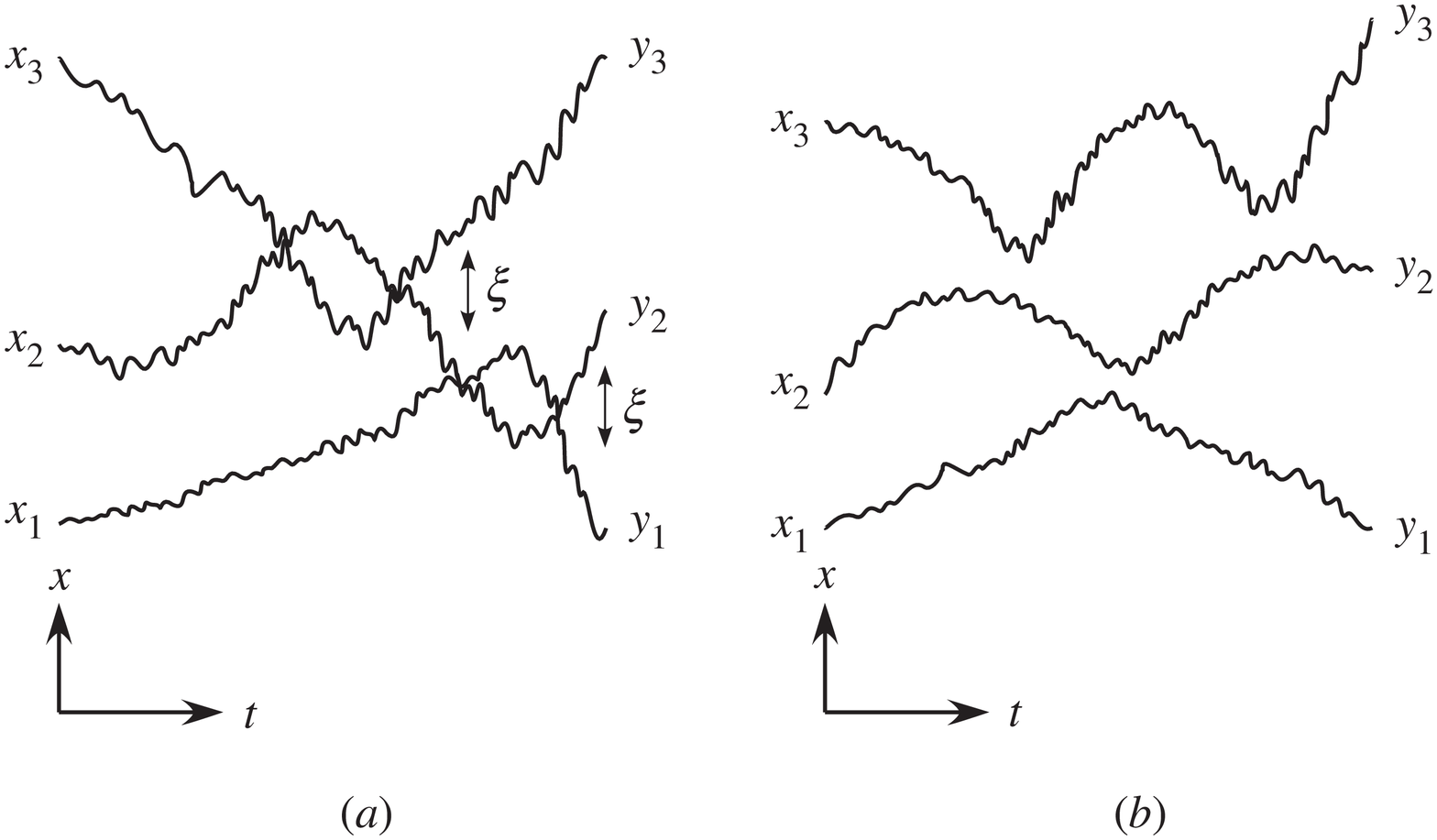}
\caption{
(1) An illustration of three paths of surviving Brownian 
particles in the mutually killing Brownian motions.
The effect of killing is not negligible when the distance between
the nearest neighbor pair of particles
is less than $\xi$;
$B_{j+1}(t) - B_j(t) < \xi$.
The ordering of particle positions
can be changed, but the risk of pair annihilation
becomes large exponentially, when
$B_j(t) - B_{j+1}(t) > \xi$.
(b) In the limit $\xi \to 0$, the present system
converges to the vicious Brownian motion,
in which the process survives if and only if paths are not
intersecting.
\label{fig:Fig2}}
\end{figure}

For $N \in \{2,3, \dots\}$, we consider the $N$-particle system
of mutually killing Brownian motions, in which the probability
that all $N$ particles survive
decays in time following (\ref{eqn:PN1}) with (\ref{eqn:VN1})
depending on realization of paths. For $\x, \y \in \R^N, t \geq 0$,
the transition probability density from the state $\x$ to $\y$ 
during time interval $t$ is then given by
\begin{eqnarray}
&& Q_N(t, \y|\x) = 
\bE \left[ 
\exp \left\{ -\int_0^t V_N(\B(s)) ds \right\} 
\1(\B(0)=\y, \B(t)=\x) \right]
\nonumber\\
&& \quad = \bE \left[
\exp \left\{ - \frac{1}{\xi^2} \sum_{j=1}^{N-1} \int_0^t
e^{-(B_{j+1}(s)-B_j(s))/\xi} ds \right\} 
\1(\B(0)=\y, \B(t)=\x) \right].
\label{eqn:QN1}
\end{eqnarray}
It is the Feynman-Kac formula for the solution of the
diffusion equation (the backward Kolmogorov equation)
\begin{equation}
\frac{\partial u(t,\x)}{\partial t}
=\frac{1}{2} \Delta u(t, \x)
-V_N(\x) u(t, \x)
\label{eqn:diff1}
\end{equation}
with $\Delta=\sum_{j=1}^{N} \partial^2/\partial x_j^2$
under the initial condition
\begin{equation}
u(0,\x)=\delta(\x-\y) \equiv
\prod_{j=1}^{N} \delta(x_j-y_j).
\label{eqn:diff2}
\end{equation}

\vskip 0.3cm
\noindent{\bf Remark 2.} 
If we regard (\ref{eqn:diff1}) as a diffusion equation
describing an $N$-dimensional Brownian motion
in $\R^N$ with a killing term $-V_N(\x)$,
(\ref{eqn:QN1}) gives a transition probability density
for a particle assumed to be at a position
$\x \in \R^N$ such that it survives during time $t$
and it arrives at a position $\y \in \R^N$
after the time duration $t$
(the Feynman-Kac formula, see, for instance, \cite{KS91}).
In the present paper, on the other hand, 
we would like to consider an $N$-particle system of 
one-dimensional Brownian motions, such that
the $j$-th and $(j+1)$-th particles will be
pair annihilated with high probability
if $B_{j+1}(t)-B_{j}(t) < \xi, 1 \leq j \leq N-1$,
and then the probability that all $N$ particles
survive avoiding from any mutual killing
decays in time following (\ref{eqn:PN1}) with (\ref{eqn:VN1}).
In order to discuss processes, in which 
mutual killing of particles actually occurs and
total number of particles decreases in time, 
we have to specify the stochastic rule
to determine which pair of particles
is annihilated, when the survival probability
conditional on a path,
$\exp\{-\int_0^t V_N(\B(s)) ds \}$, becomes small.
Here we are interested in, however, the situation
that mutual killing of particles does not occur at all
and all $N$ particles survive, following the
notion of vicious walker models 
\cite{Fis84,KT02,CK03,Kat11}.
The Feynman-Kac formula (\ref{eqn:QN1}) gives the
transition probability density between
$N$-particle configurations $\x, \y \in \R^N$.
Figure 2(a) illustrates three paths of surviving
particles, in which change of ordering
of particle positions occurs within
the spatial scale $\simeq \xi$.
\vskip 0.3cm

By the fact that the Whittaker function $\psi_{\k}^{(N)}(\x)$
solves (\ref{eqn:H2}) with (\ref{eqn:lambda1}),
we can see that
\begin{equation}
\exp \left( -\frac{t}{2} |\k|^2 \right) 
\psi_{i \xi \k}^{(N)}(\x/\xi),
\quad 
|\k|^2 \equiv \sum_{j=1}^{N} k_j^2,
\quad \k \in \R^N, 
\label{eqn:sol1}
\end{equation}
solves the diffusion equation (\ref{eqn:diff1}).

The density function of the Sklyanin measure \cite{Skl85}
is defined by
\begin{equation}
s_N(\k)=\frac{1}{(2 \pi)^N N!} 
\prod_{1 \leq j < \ell \leq N}
|\Gamma(i(k_{\ell}-k_{j}))|^{-2}
\label{eqn:Skl1}
\end{equation}
for $\k \in \C^N$.
By Euler's reflection formula
$\Gamma(z) \Gamma(1-z)=\pi/\sin(\pi z)$,
if $\k \in \R^N$,
$$
\prod_{1 \leq j < \ell \leq N}
|\Gamma(i(k_{\ell}-k_j))|^{-2}
=\prod_{1 \leq j < \ell \leq N}
\left\{ (k_{\ell}-k_j)
\frac{\sinh \pi(k_{\ell}-k_j)}{\pi} \right\}.
$$
Since the Whittaker functions satisfy the completeness
relation with respect to the Sklyanin measure
$s_N(\k) d\k \equiv s_N(\k) \prod_{j=1}^{N} dk_j$
\cite{STS94,KL99,KL01,GLO08},
\begin{equation}
\int_{\R^N} \psi^{(N)}_{i \k}(\x)
\psi^{(N)}_{-i \k}(\y) s_N(\k) d \k
=\delta(\x-\y), \quad \x, \y \in \R^N, 
\label{eqn:complete1}
\end{equation}
the integral of the solution (\ref{eqn:sol1}) 
multiplied by $\psi^{(N)}_{-i \xi \k}(\y/\xi)$ 
on this measure, 
\begin{equation}
Q_N(t,\y|\x)=\int_{\R^N}
e^{-t|\k|^2/2} \psi_{i \xi \k}^{(N)}(\x/\xi)
\psi_{-i \xi \k}^{(N)}(\y/\xi) s_N(\xi \k) d \k
\label{eqn:QN2}
\end{equation}
satisfies the initial condition (\ref{eqn:diff2}).
That is, (\ref{eqn:QN2}) is an integral representation
of the transition probability density (\ref{eqn:QN1}).
It should be noted that the orthogonality relation
\cite{KL99,GLO08}
\begin{equation}
\int_{\R^N} \psi_{-i \k}^{(N)}(\x) \psi_{i \k'}^{(N)}(\x) d \x
=\frac{1}{s_N(\k) N!} \sum_{\sigma \in \cS_N}
\delta(\k-\sigma(\k')),
\quad \k,\k' \in \R^N,
\label{eqn:ortho1}
\end{equation}
is established, 
where $\cS_N$ is a set of all permutation of $N$
indices and $\sigma(\k') \equiv (k'_{\sigma(1)}, \dots, k'_{\sigma(N)})$.
It guarantees the Chapman-Kolmogorov equation
for the transition probability density
\begin{equation}
\int_{\R^N} Q_N(t_2, \z|\y) Q_N(t_1, \y|\x) d \y
=Q_N(t_1+t_2, \z|\x),
\quad \x, \z \in \R^N
\label{eqn:CK1}
\end{equation}
for $0 \leq t_1, t_2 < \infty$, which should hold,
since the present process is Markovian.

\vskip 0.3cm
\noindent{\bf Remark 3.} \,
The completeness relation (\ref{eqn:complete1})
of the Whittaker functions is a multivariate
extension of the relation for Macdonald's functions,
\begin{equation}
\frac{1}{\pi^2} \int_{\R} K_{i k}(e^{-x})
K_{-i k}(e^{-y}) k \sinh(\pi k) d k
=\delta(x-y), \quad x, y \in \R, 
\label{eqn:complete2}
\end{equation}
which will define the Kontorovich-Lebedev transform 
\cite{Leb65}.
We can use (\ref{eqn:complete2})
to confirm that (\ref{eqn:Q3}) satisfies the 
initial condition $Q(0, y|x)=\delta(x-y)$.
\vskip 0.3cm

\subsection{Asymptotics of transition probability density
and survival probability}

By using the integral representation (\ref{eqn:WInt1}),
we can prove the following asymptotics of the
transition probability density.

\begin{lem}
\label{thm:AsymQ}
For $\x, \y \in \R^N, |\x| < \infty, t >0$, 
as $|\x|/\sqrt{t} \to 0$, 
\begin{eqnarray}
&& Q_N(t, \y|\x) = \xi^{N(N-1)}
\frac{1}{(2\pi)^N N!} \left(\frac{2}{t}\right)^{N^2/2}
\psi_0^{(N)}(\x/\xi) \psi_0^{(N)}(\y/\xi)
e^{-|\y|^2/2t}
\nonumber\\
&& \quad \times \int_{\R^N} e^{-|\vmu|^2}
\prod_{1 \leq j< k \leq N}
\left[ (\mu_k-\mu_j)^2+\frac{1}{2t}(y_k-y_j)^2 \right] d \vmu
\times \left\{ 1+ {\cal O}
\left( \frac{|\x|}{\sqrt{t}} \right) \right\}.
\label{eqn:AsymQ}
\end{eqnarray}
\end{lem}
\vskip 0.5cm
\noindent
On the other hand, by using Givental's integral
representation (\ref{eqn:WInt2}), we have the
following \cite{OCo09}. Let $\W_N$ be the Weyl chamber of type A$_{N-1}$,
$$
\W_N=\{\x=(x_1, x_2, \dots, x_N) \in \R^N :
x_1 < x_2 < \cdots < x_N \}
$$
and set
\begin{equation}
h_N(\x)=\prod_{1 \leq j < k \leq N}(x_k-x_j).
\label{eqn:hN1}
\end{equation}

\begin{lem}
\label{thm:Asympsi0}
For $\x \in \R^N$,
\begin{equation}
\lim_{\beta \to \infty} \beta^{-N(N-1)/2}
\psi_0^{(N)}(\beta \x)
= \frac{h_N(\x)}{\prod_{j=1}^{N-1} j !} \1(\x \in \W_N).
\label{eqn:Asympsi0}
\end{equation}
\end{lem}
\vskip 0.3cm

The probability that all $N$ particles survive
up to time $t$ in this 
system of the mutually killing Brownian motions
is given by
\begin{equation}
\cN_N(t,\x) = \int_{\R^N} Q_N(t,\y|\x) d\y
\label{eqn:NN1}
\end{equation}
for each initial state $\x \in \R^N$.
We call it simply the {\it survival probability}
in this paper.
The following long-term asymptotics of
the survival probability is obtained.

\begin{prop}
\label{thm:AsymN}
Let 
\begin{equation}
C_N = \xi^{N(N-1)/2} 
\frac{2^{3N(N-1)/4}}{\pi^N \prod_{j=1}^{N-1} j!} A_N
\label{eqn:CN}
\end{equation}
with
\begin{equation}
A_N = 
\int_{\W_N} d \veta \int_{\W_N} d \vmu \,
e^{-(|\veta|^2+|\vmu|^2)}  
\prod_{1 \leq j < k \leq N}
\Big[ (\eta_k-\eta_j)
\{(\mu_k-\mu_j)^2+(\eta_k-\eta_j)^2\} \Big].
\label{eqn:AN}
\end{equation}
Then for $\x \in \R^N$,
\begin{equation}
\lim_{t \to \infty} t^{N(N-1)/4}
\cN_N(t, \x)=C_N \psi_0^{(N)}(\x/\xi).
\label{eqn:AsymN}
\end{equation}
\end{prop}
\vskip 0.3cm

\subsection{O'Connell process as a system of 
mutually killing Brownian motions conditioned 
that all particles survive forever}

Let $0 \leq T < \infty$ and we consider the $N$-particle
system of the mutually killing Brownian motions
conditioned that all particles survive up to time $T$.
For $0 \leq s \leq t \leq T$, the transition probability density
from the state $(s, \x)$ to $(t, \y), \x, \y \in \R^N$
of this conditional process is given by
\begin{equation}
P_N^{T}(s, \x; t, \y)
=\frac{\cN_N(T-t,\y)}{\cN_N(T-s,\x)}
Q_N(t-s, \y|\x),
\quad 0 \leq s \leq t \leq T.
\label{eqn:PNT1}
\end{equation}
The fact that (\ref{eqn:PNT1}) depends
not only $t-s$ but also $T-t$ and $T-t$ implies that
this conditional process is inhomogeneous in time.

Proposition \ref{thm:AsymN} enables us to take the
limit $T \to \infty$ of (\ref{eqn:PNT1}).
As a result we obtain the temporally homogeneous process
with the transition probability density
\begin{eqnarray}
P_N(t, \y|\x) &\equiv& \lim_{T \to \infty}
P_N^T(0,\x;t,\y)
\nonumber\\
&=& \frac{\psi_0^{(N)}(\y/\xi)}{\psi_0^{(N)}(\x/\xi)}
Q_N(t,\y|\x), \quad
\x,\y \in \R^N, t \geq 0.
\label{eqn:PN2}
\end{eqnarray}

It is easy to confirm the following \cite{Kat11}.

\vskip 0.3cm
\begin{prop}
\label{thm:diffeq}
The transition probability density (\ref{eqn:PN2}) of the
$N$-particle system of the mutually killing Brownian motions
conditioned that all particles survive forever solves the following 
diffusion equation, 
\begin{eqnarray}
&& \frac{\partial}{\partial t} u(t, \x)
= \frac{1}{2} \Delta u(t, \x)
+\nabla \log \psi_0^{(N)}(\x/\xi) \cdot \nabla u(t, \x)
\nonumber\\
&& = \frac{1}{2} \sum_{j=1}^{N}
\frac{\partial^2}{\partial x_j^2} u(t, \x)
+\sum_{1 \leq j \leq N}
\frac{\partial \log \psi_0^{(N)}(\x/\xi)}{\partial x_j}
\frac{\partial}{\partial x_j} u(t, \x),
\quad \x \in \R^N, t \geq 0, 
\label{eqn:OConnell1}
\end{eqnarray}
under the initial condition
$u(0, \x)=\delta(\x-\y), \y \in \R^N$.
\end{prop}
\vskip 0.3cm

Recently O'Connell \cite{OCo09} introduced a diffusion process
in $\R^N$ with the infinitesimal generator of the process
\begin{eqnarray}
\cL^{\vnu}_N &=&
- (\psi^{(N)}_{\vnu})^{-1} \Big(\cH_N +|\vnu|^2 \Big) \psi^{(N)}_{\vnu}
\nonumber\\
&=& \frac{1}{2} \Delta
+\nabla \log \psi^{(N)}_{\vnu}(\x) \cdot \nabla
\label{eqn:OConnell3}
\end{eqnarray}
with drift $\vnu \in \R^N$.
Then we can state the following.

\vskip 0.3cm
\begin{thm}
\label{thm:OConnell}
The O'Connell process without drift, $\vnu=0$, is
realized as the system of mutually killing Brownian motions
conditioned that all particles survive forever.
\end{thm}
\vskip 0.5cm

\noindent
In an earlier paper \cite{Kat11},
we considered the system of mutually killing
Brownian motions with drifts
$\vnu=(\nu_1, \dots, \nu_N) \in \R^N$
conditioned that all particles survive forever
and took the limit $\vnu \to 0$
to realize the O'Connell process
without drift $\vnu=0$.
In the present paper we have shown
that an introduction of drift $\vnu$
is not necessary to deriving (\ref{eqn:PN2})
and concluding Theorem \ref{thm:OConnell}.

\subsection{The $\xi \to 0$ limit}

It is known that the functions $\{m^{(N)}(\x, \k)\}_{\k \in \R^N}$ called
the fundamental Whittaker functions are defined
so that the (class-one) Whittaker functions
discussed so far are expressed by the 
following alternating sum of them \cite{KL99,BO11},
\begin{equation}
\psi^{(N)}_{i \k}(\x) = \prod_{1 \leq j < \ell \leq N}
\frac{\pi}{i \sinh \pi (k_{\ell}-k_j)}
\sum_{\sigma \in \cS_N} {\rm sgn}(\sigma)
m^{(N)}(\x, \sigma(\k)),
\label{eqn:alt1}
\end{equation}
where $\sigma(\k)=(k_{\sigma(1)}, \dots, k_{\sigma(N)})$
for each permutation $\sigma \in \cS_N$.
The functions $\{m^{(N)}(\x, \k)\}_{\k \in \R^N}$ 
are normalized as
\begin{equation}
\lim_{\xi \to 0} m^{(N)}(\x/\xi, \xi \k)
=e^{i \k \cdot \x} \quad
\mbox{for $\x \in \W_N$}.
\label{eqn:mN1}
\end{equation}
Note that (\ref{eqn:alt1}) is a multivariate extension
of the definition (\ref{eqn:Knu1}) of $K_{\nu}(z)$, which is given by
an alternating sum of $I_{\nu}(z)$ and $I_{-\nu}(z)$.
Therefore, when $\x, \y \in \W_N$, 
we can see that (\ref{eqn:QN2}) with (\ref{eqn:Skl1})
gives 
\begin{eqnarray}
q_N(t, \y|\x) &\equiv& \lim_{\xi \to 0} Q_N(t, \y|\x)
\nonumber\\
&=& \frac{1}{(2 \pi)^N N!} \sum_{\sigma \in \cS_N}
\sum_{\sigma' \in \cS_N} {\rm sgn}(\sigma){\rm sgn}(\sigma')
\int_{\R^N} e^{-t|\k|^2/2+i(\sigma(\k)\cdot \x-\sigma'(\k)\cdot \y)} d \k
\nonumber\\
&=& \frac{1}{(2 \pi)^N} \sum_{\rho \in \cS_N} {\rm sgn}(\rho)
\int_{\R^N} e^{-t|\k|^2/2+i \k \cdot (\x-\rho(\y))} d \k
\nonumber\\
&=& \frac{1}{(2 \pi)^N} \sum_{\rho \in \cS_N}
{\rm sgn}(\rho) 
\left( \frac{2 \pi}{t} \right)^{N/2} e^{-|\x-\rho(\y)|^2/2t}
\nonumber\\
&=& \det_{1 \leq j, k \leq N}
\left[ \frac{1}{\sqrt{2 \pi t}}
e^{-(x_j-y_k)^2/2t} \right].
\label{eqn:KM1}
\end{eqnarray}
This is the Karlin-McGregor determinant 
giving a total mass of nonintersecting paths
of Brownian motions from $\x \in \W_N$
to $\y \in \W_N$ during time $t$ \cite{KM59,KT_Sugaku_11}.
As a matter of fact, when $\x, \y \in \W_N$, 
the limit $\xi \to 0$ of (\ref{eqn:QN1})
should give the transition probability density of the
$N$-dimensional absorbing Brownian motion in the Weyl chamber $\W_N$
(the vicious Brownian motion \cite{Kat11}).
See Fig.\ref{fig:Fig2}(b).
Lemma \ref{thm:Asympsi0} implies that, if $\x, \y \in \W_N$,
\begin{eqnarray}
p_N(t,\y|\x) &\equiv& \lim_{\xi \to 0} P_N(t, \y|\x)
\nonumber\\
&=& \frac{h_N(\y)}{h_N(\x)} q_N(t, \y|\x),
\quad t \geq 0,
\label{eqn:pN1}
\end{eqnarray}
which is an $h$-transform of $q_N$
by the Vandermonde determinant (\ref{eqn:hN1}).
It shows that the limit $\xi \to 0$ of the
O'Connell process is identified with the Dyson model
(Dyson's Brownian motion model with $\beta=2$ \cite{Dys62}),
which is the $N$-particle system of Brownian motions
conditioned never to collide with each other
(the noncolliding Brownian motion) \cite{KT02,KT10,KT_Sugaku_11}.

\vskip 0.3cm
\noindent{\bf Remark 4.} \
The survival probability up to time $T < \infty$
of the one-dimensional absorbing Brownian motion
with an absorbing wall at the origin is given by
\begin{equation}
\cN^0(T,x)=\int_{-\infty}^{\infty} q(T,y|x) dy,
\label{eqn:NZ0}
\end{equation}
where the transition probability density $q$ is
given by (\ref{eqn:q1}). The long-term behavior
is readily obtained as
\begin{equation}
\cN^0(T,x) \simeq c^0 x T^{-\phi}, \quad
x>0, \quad T \to \infty
\label{eqn:NZ1}
\end{equation}
with
$c^0=\sqrt{2/\pi},
\phi=1/2$.
Therefore, the exponent
$\phi$ is the same with that found in the
asymptotics (\ref{eqn:Nlimit}) of the survival probability
of the killing Brownian motion discussed in Sect.1;
\begin{equation}
\cN(T,x) \simeq 3 \sqrt{\frac{2}{\pi}} \xi
K_0(e^{-x/\xi}) T^{-\phi}, \quad
x \in \R, \quad T \to \infty.
\label{eqn:NZ3}
\end{equation}
If we take the limit $\xi \to 0$,
(\ref{eqn:KlimB}) gives
\begin{equation}
\lim_{\xi \to 0} \cN(T,x)
\simeq c x T^{-\phi}, \quad
x > 0, \quad T \to \infty
\label{eqn:NZ4}
\end{equation}
with
$c=3 \sqrt{2/\pi} = 3 c^0$.
For $N \geq 2$, Proposition \ref{thm:AsymN} implies
\begin{equation}
\cN_N(T,\x) \simeq C_N \psi^{(N)}_0(\x/\xi)
T^{-\phi_N}, \quad
\x \in \R^N, \quad T \to \infty
\label{eqn:NNb1}
\end{equation}
with $C_N$ given by (\ref{eqn:CN}) and with the
{\it survival probability exponent}
\begin{equation}
\phi_N=\frac{1}{4} N(N-1)
\label{eqn:NNb2}
\end{equation}
for the present $N$-particle system of 
mutually killing Brownian motion.
This exponent is the same as that
for the vicious Brownian motion
(the absorbing Brownian motion in
$\W_N$) \cite{Fis84,KT02};
\begin{eqnarray}
\cN^0_N(T,\x) &\equiv& \int_{\W^N} q_N(T,\y|\x) d \y
\nonumber\\
&\simeq& c_N^0 h_N(\x) T^{-\phi_N},
\quad \x \in \W_N, \quad T \to \infty
\label{eqn:NNb3}
\end{eqnarray}
with 
\begin{equation}
c_N^0=\frac{2^{N/2} \prod_{j=1}^{N} \Gamma(j/2)}
{\pi^{N/2} \prod_{j=1}^{N-1} j!}.
\label{eqn:NNb4}
\end{equation}
Lemma \ref{thm:Asympsi0} gives then
\begin{equation}
\lim_{\xi \to 0} \cN_N(T,\x)
\simeq c_N h_N(\x) T^{-\phi_N},
\quad \x \in \W_N, \quad T \to \infty,
\label{eqn:NNb5}
\end{equation}
where
\begin{equation}
c_N=\frac{2^{3N(N-1)/4}}{\pi^N 
\prod_{j=1}^{N-1}(j!)^2} A_N
\label{eqn:NNb6}
\end{equation}
with (\ref{eqn:AN}).
The above results $c \not= c^0$ and $c_N \not= c_N^0,
N \geq 2$ are consequences of the facts that
for the transition probability densities
$Q(t,y|x)$ and $Q_N(t,\y|\x), N \geq 2$,
the long-term limit $T \to \infty$ and 
the limit $\xi \to 0$ 
are noncommutable.

\SSC{Proofs of Lemmas and Proposition}
\subsection{Proof of Lemma \ref{thm:AsymQ}}

By the integral representation (\ref{eqn:WInt1}), we put
$\psi^{(N)}_{i \xi \k}(\x/\xi)=\exp(i \k \cdot \x)
\widehat{\psi}^{(N)}_{i \xi \k}(\x/\xi)$.
Then by changing the integral variables as
$\mu_{j}=\sqrt{t/2} k_j-i(x_j-y_j)/\sqrt{2t}, 1 \leq j \leq N$,
(\ref{eqn:QN2}) with (\ref{eqn:Skl1}) is written as
\begin{eqnarray}
&& Q_N(t, \y|\x) = \frac{1}{(2\pi)^N N!}
\left( \frac{2}{t} \right)^{N/2} e^{-|\x-\y|^2/2t}
\nonumber\\
&& \quad \times \int_{\R^N}
e^{-|\vmu|^2} \widehat{\psi}^{(N)}_{i \valpha(\vmu)}(\x/\xi)
\widehat{\psi}^{(N)}_{-i \valpha(\vmu)}(\y/\xi)
\prod_{1 \leq j < k \leq N}
\Big|\Gamma(i (\alpha_k(\mu_k)-\alpha_j(\mu_j))) \Big|^{-2} d \vmu,
\nonumber
\end{eqnarray}
where 
$$
\alpha_j(\mu_j)=\xi \sqrt{\frac{2}{t}} \left\{
\mu_j+\frac{i}{\sqrt{2t}} (x_j-y_j) \right\},
\quad 1 \leq j \leq N,
$$
and $\valpha(\vmu)=(\alpha_1(\mu_1), \dots,
\alpha_N(\mu_N))$.
For $|\x| < \infty$, in the limit
$|\x|/\sqrt{t} \to 0$, 
$\alpha_j(\mu_j) \simeq \xi \sqrt{2/t}
\{ \mu_j-i y_j/\sqrt{2t} \} \to 0, 1 \leq j \leq N$,
and thus $\widehat{\psi}^{(N)}_{i \valpha(\vmu)}(\x/\xi)
\to \widehat{\psi}^{(N)}_0(\x/\xi) \equiv
\psi^{(N)}_{0}(\x/\xi)$.
Since $z \Gamma(z) \to 1$ as $z \to 0$,
$$
\Big| \Gamma(i(\alpha_k(\mu_k)-\alpha_j(\mu_j))) \Big|^{-2}
\simeq \xi^2 \frac{2}{t} \left\{
(\mu_k-\mu_j)^2+\frac{1}{2t} (y_k-y_j)^2 \right\}
$$
in this limit. Then (\ref{eqn:AsymQ}) is obtained. \qed

\subsection{Proof of Lemma \ref{thm:Asympsi0}}

The following argument is found in Section 6 of O'Connell \cite{OCo09}.
For $N \geq 2$, the {\it Gelfand-Tsetlin pattern}
({\it complete branching}) is defined as \cite{Stanley99}
\begin{eqnarray}
\GT_N &=& 
\Big\{ (\x^1, \x^2, \dots, \x^N) \in \overline{\W}_1 \times
\overline{\W}_2 \times \cdots \times \overline{\W}_N :
\nonumber\\
&& \qquad \quad
x^{j+1}_{k} \leq x^{j}_{k} \leq x^{j+1}_{k+1},
1 \leq k \leq j \leq N-1 \Big\},
\label{eqn:GT}
\end{eqnarray}
where $\overline{\W}_j=\{\x=(x_1, \dots, x_j) \in \R^j:
x_1 \leq \cdots \leq x_j \}, 1 \leq j \leq N$.
By the Givental integral representation (\ref{eqn:WInt2}),
\begin{eqnarray}
\psi_{0}^{(N)}(\beta \x)
&=&\int_{\Gamma_N(\beta \x)} 
\exp \left[ 
- \sum_{j=1}^{N-1} \sum_{k=1}^j
\Big\{ e^{-(T_{j,k}-T_{j+1,k})}
+e^{-(T_{j+1, k+1}-T_{j,k})} \Big\} \right] d \T
\nonumber\\
&=& \beta^{N(N-1)/2} \int_{\Gamma_N(\x)} 
\exp \left[ 
- \sum_{j=1}^{N-1} \sum_{k=1}^j
\Big\{ e^{-\beta({T'}_{j,k}-{T'}_{j+1,k})}
+e^{-\beta({T'}_{j+1, k+1}-{T'}_{j,k})} \Big\} \right] d \T',
\nonumber
\end{eqnarray}
where we set 
$\T=\beta \T'=(\beta {T'}_{j,k}: 1 \leq k \leq j \leq N)$.
Let $\t^j=({T'}_{j,1}, \dots, {T'}_{j,j}), 1 \leq j \leq N$.
As $\beta \to \infty$, if
\begin{equation}
(\t^1, \t^2, \dots, \t^N) \in \GT_N \quad
\mbox{and} \quad \t^N=\x \in \W_N,
\label{eqn:GT2}
\end{equation}
then the integrand converges to 1, and
otherwise it becomes zero.
Then the integral is the volume of $\GT_N$ under the
condition $\x^N=\x \in \W_N$ (the GT-polytope),
which is given by 
$h_N(\x)/\prod_{j=1}^{N-1} j!$.
Then the proof is completed. \qed

\subsection{Proof of Proposition \ref{thm:AsymN}}

By definition (\ref{eqn:NN1}), Lemma \ref{thm:AsymQ} gives
\begin{eqnarray}
&& \cN_N(t, \x) = \xi^{N(N-1)} \frac{1}{(2 \pi)^N N!}
\left( \frac{2}{t} \right)^{N^2/2}
\psi_0^{(N)}(\x/\xi)
\int_{\R^N} d \y \, e^{-|\y|^2/2t} \psi^{(N)}_0(\y/\xi)
\nonumber\\
&& \qquad \times
\int_{\R^N} d \vmu \, e^{-|\vmu|^2}
\prod_{1 \leq j < k \leq N}
\left[ (\mu_k-\mu_k)^2+\frac{1}{2t}(y_k-y_j)^2 \right]
\times \left\{ 1+{\cal O} \left( \frac{|\x|}{\sqrt{t}} \right) \right\}
\nonumber
\end{eqnarray}
in $|\x|/\sqrt{t} \to 0$.
If we put $y_j/\sqrt{2t} = \eta_j, 1 \leq j \leq N$,
it is written as
\begin{eqnarray}
\cN_N(t, \x) &=& \xi^{N(N-1)} \frac{1}{\pi^N N!}
\left( \frac{2}{t} \right)^{N(N-1)/2} \psi_0^{(N)}(\x/\xi)
\nonumber\\
&& \times
\int_{\R^N} d \veta 
\int_{\R^N} d \vmu \,
e^{-(|\veta|^2+|\vmu|^2)}
\psi^{(N)}_{0} \left( \frac{\sqrt{2t}}{\xi} \veta \right)
\nonumber\\
&& \times \prod_{1 \leq j < k \leq N}
\Big[ (\mu_k-\mu_j)^2+(\eta_k-\eta_j)^2 \Big]
\times \left\{ 1+ {\cal O}
\left( \frac{|\x|}{\sqrt{t}} \right) \right\},
\nonumber
\end{eqnarray}
in $|\x|/\sqrt{t} \to 0$.
By Lemma \ref{thm:Asympsi0}, 
$$
\psi^{(N)}_0 \left( \frac{\sqrt{2t}}{\xi} \veta \right)
\simeq \left( \frac{\sqrt{2t}}{\xi} \right)^{N(N-1)/2}
\frac{h_N(\veta)}{\prod_{j=1}^{N-1} j !} \1 ( \veta \in \W_N),
$$
in $t \to \infty$, and then
the proposition is proved. \qed

\clearpage
\appendix
\begin{LARGE}
{\bf Appendix}
\end{LARGE}
\SSC{One-dimensional killing Brownian motion}

By (\ref{eqn:KInt1}), we can put $K_{i \xi k}(e^{-x/\xi}) 
=e^{i k x} \widehat{K}_{i \xi k}(e^{-x/\xi})$
with $\widehat{K}_{\nu}(z)=(2^{\nu} \Gamma(\nu+1/2)/\sqrt{\pi})$
$\int_0^{\infty} du \cos(zu)/(1+u^2)^{\nu+1/2}$.
Changing the integral variable in (\ref{eqn:Q3}) by
$k \sqrt{t/2} -i(x-y)/\sqrt{2t}=\mu$, 
the transition probability density can be written as
\begin{eqnarray}
&& Q(t,y|x)= \frac{1}{\pi} \sqrt{\frac{2}{t}} e^{-(x-y)^2/2t}
\nonumber\\
&& \quad \quad \times
\int_{-\infty}^{\infty} e^{-\mu^2} 
\widehat{K}_{i \alpha(\mu)}(e^{-x/\xi})
\widehat{K}_{- i \alpha(\mu)}(e^{-y/\xi})
|\Gamma(i \alpha(\mu))|^{-2} d \mu
\nonumber
\end{eqnarray}
with
$\alpha(\mu) =\xi \sqrt{2/t} \{\mu+i(x-y)/\sqrt{2t}\}$.
For $|x| < \infty$, as $x/\sqrt{t} \to 0$,
$\alpha(\mu) \simeq \xi \sqrt{2/t}(\mu-iy/\sqrt{2t})
\to 0$, and
$\Gamma(i \alpha(\mu))
\simeq \{i \xi \sqrt{2/t} (\mu-iy/\sqrt{2t}) \}^{-1}$.
Then 
\begin{eqnarray}
&& Q(t, y|x) =
\frac{\xi^2}{\pi} \left( \frac{2}{t} \right)^{3/2}
e^{-y^2/2t} K_0(e^{-x/\xi}) K_0(e^{-y/\xi})
\nonumber\\
&& \qquad \qquad \times
\int_{-\infty}^{\infty} e^{-\mu^2} 
\left( \mu^2 +\frac{y^2}{2t} \right) d \mu
\times \left\{ 1+ {\cal O}
\left( \frac{x}{\sqrt{t}} \right) \right\}
\nonumber\\
&& \qquad = \xi^2 \sqrt{\frac{2}{\pi}}
t^{-3/2} e^{-y^2/2t}
\left( 1+\frac{y^2}{t} \right)
K_0(e^{-x/\xi}) K_0(e^{-y/\xi})
\times 
 \left\{ 1+ {\cal O} \left( \frac{x}{\sqrt{t}} \right) \right\}
\nonumber
\end{eqnarray}
in $x/\sqrt{t} \to 0$, since $\widehat{K}_0=K_0$ by definition.
It gives through (\ref{eqn:N1})
\begin{equation}
\cN(t,x) = \xi^2 \sqrt{\frac{2}{\pi}}
t^{-3/2} K_0(e^{-x/\xi}) I_t \times 
\left\{ 1+ {\cal O} \left( \frac{x}{\sqrt{t}} \right) \right\}
\label{eqn:Nasym}
\end{equation}
with
\begin{eqnarray}
I_t &=& \int_{-\infty}^{\infty} e^{-y^2/2t}
\left( 1 + \frac{y^2}{t} \right)
K_0(e^{-y/\xi}) dy
\nonumber\\
&=& \sqrt{2t}
\int_{-\infty}^{\infty} e^{-u^2}
(1+2u^2) K_0(e^{-\sqrt{2t} u/\xi}) du.
\nonumber
\end{eqnarray}

By (\ref{eqn:KInt2}), for $\beta >0$,
$K_0(e^{-\beta x})$ is written as
\begin{eqnarray}
K_0(e^{-\beta x}) &=&
\frac{1}{2} \int_0^{\infty}
s^{-1} \exp \left\{ -\frac{1}{2} e^{-\beta x}
\left( s + \frac{1}{s} \right) \right\} ds
\nonumber\\
&=& \frac{\beta}{2} \int_{-\infty}^{\infty}
\exp \left\{ -\frac{1}{2}
\Big( e^{-\beta(x-v)}+e^{-\beta(x+v)} \Big) \right\} dv,
\nonumber
\end{eqnarray}
where we have set $s=e^{\beta v}$.
In the limit $\beta \to \infty$,
$e^{-\beta(x-v)}+e^{-\beta(x+v)} =0$, 
if $v < x$ and $v>-x$, 
and $e^{-\beta(x-v)}+e^{-\beta(x+v)} = \infty$, 
otherwise.
Therefore
\begin{equation}
\lim_{\beta \to \infty} \beta^{-1}
K_0(e^{-\beta x})
= \frac{1}{2} \int_{-x}^{x} dv \1(x>0)
= x \, \1(x>0).
\label{eqn:K0limA}
\end{equation}

Then $K_0(e^{-\sqrt{2t}u/\xi}) \simeq
\sqrt{2t}u/\xi \1(u > 0)$ in $t \to \infty$, and 
$$
I_t = \frac{2t}{\xi} \int_{0}^{\infty}
e^{-u^2} (1+2u^2) u \, du \times
\{1+o(1) \}
=\frac{3t}{\xi} \times \{1+o(1)\}
\quad \mbox{in} \quad t \to \infty.
$$
Put this result into (\ref{eqn:Nasym}),
then we obtain the asymptotics (\ref{eqn:Nlimit}).

\vskip 0.5cm
\noindent{\bf Acknowledgements} \quad
The present author would like to thank
T. Sasamoto and T. Imamura 
for useful discussion on the present work.
A part of the present work was done
during the participation of the present author 
in \'Ecole de Physique des Houches on
`` Vicious Walkers and Random Matrices"
(May 16-27, 2011).
The author thanks G. Schehr, C. Donati-Martin, and
S. P\'ech\'e 
for well-organization of the school.
This work is supported in part by
the Grant-in-Aid for Scientific Research (C)
(No.21540397) of Japan Society for
the Promotion of Science.



\begin{thebibliography}{99}

\bibitem{BO11}
Baudoin, F., O'Connell, N.:
Exponential functionals of Brownian motion
and class-one Whittaker functions.
Ann. Inst. H. Poincar\'e, 
B {\bf 47}, 1096-1120 (2011)

\bibitem{BC11}
Borodin, A., Corwin, I.:
Macdonald processes.
{\sf arXiv:math.PR/1111.4408}

\bibitem{CK03}
Cardy, J., Katori, M.:
Families of vicious walkers.
J. Phys. A {\bf 36}, 609-629 (2003)

\bibitem{COSZ11}
Corwin, I., O'Connell, N., Sepp\"al\"ainen, T.,
Zygouras, N.:
Tropical combinatorics and Whittaker functions.
{\sf arXiv:math.PR/1110.3489}

\bibitem{Dys62}
Dyson, F. J.: 
A Brownian-motion model for the eigenvalues of a random
matrix.
J. Math. Phys. {\bf 3}, 1191-1198 (1962)

\bibitem{Fis84}
Fisher, M. E.:
Walks, walls, wetting, and melting.
J. Stat. Phys. {\bf 34}, 667-729 (1984)

\bibitem{GKLO06}
Gerasimov, A., Kharchev, S., Lebedev, D., Oblezin, S.:
On a Gauss-Givental representations of quantum Toda chain
wave equation. 
Int. Math. Res. Not. 1-23 (2006)

\bibitem{GLO08}
Gerasimov, A., Lebedev, D., Oblezin, S.:
Baxter operator and Archimedean Hecke algebra.
Commun. Math. Phys. {\bf 284}, 867-896 (2008)

\bibitem{Giv97}
Givental, A.:
Stationary phase integrals, quantum Toda lattices,
flag manifolds and the mirror conjecture.
In: Topics in Singular Theory,
AMS Trans. Ser. 2, vol. 180, pp.103-115,
AMS, Rhode Island (1997)

\bibitem{GNSS11}
Gorsky, A., Nechaev, S., Santachiara, R.,
Schehr, G.:
Random ballistic growth and diffusion in 
symmetric spaces.
{\sf arXiv:math-ph/1110.3524}

\bibitem{Gra99}
Grabiner, D. J.:
Brownian motion in a Weyl chamber,
non-colliding particles, and random matrices,
Ann. Inst. Henri Poincar\'e, Probab. Stat. {\bf 35}, 177-204 (1999)

\bibitem{IS07}
Ishii, T., Stade, E.:
New formulas for Whittaker functions 
on $GL(n, \R)$.
J. Func. Anal. {\bf 244}, 289-314 (2007)

\bibitem{Jac67}
Jacquet, H.:
Fonctions de Whittaker associ\'ees aux groupes
de Chevalley.
Bull. Soc. Math. France {\bf 95}, 243-309 (1967)

\bibitem{JK03}
Joe, D., Kim, B.:
Equivariant mirrors and the Virasoro conjecture for
flag manifolds.
Int. Math. Res. Notices, no.15, 859-882 (2003)

\bibitem{KS91}
Karatzas, I., Shreve, S. E.:
Brownian Motion and Stochastic Calculus. 2nd edn.
Springer, (1991)

\bibitem{KM59}
Karlin, S., McGregor, J.:
Coincidence probabilities.
Pacific J. Math. {\bf 9}, 1141-1164 (1959)

\bibitem{Kat11}
Katori, M.:
O'Connell's process as a vicious Brownian motion.
to appear in Phys. Rev. E. 
{\sf arXiv:math-ph/1110.1845}

\bibitem{KT02}
Katori, M., Tanemura, H.:
Scaling limit of vicious walks and two-matrix model.
Phys. Rev. E {\bf 66}, 011105 (2002)

\bibitem{KT10}
Katori, M., Tanemura, H.:
Non-equilibrium dynamics of Dyson's model with
an infinite number of particles.
Commun. Math. Phys. {\bf 293}, 469-497 (2010)

\bibitem{KT_Sugaku_11}
Katori, M., Tanemura, H.:
Noncolliding processes, matrix-valued processes 
and determinantal processes.
Sugaku Expositions {\bf 24}, 263-289 (2011);
{\sf arXiv:math.PR/1005.0533}

\bibitem{KL99}
Kharchev, S., Lebedev, D.:
Integral representation for the eigenfunctions
of a quantum periodic Toda chain.
Lett. Math. Phys. {\bf 50}, 55-77 (1999)

\bibitem{KL00}
Kharchev, S., Lebedev, D.:
Eigenfunctions of $GL(N,\R)$ Toda chain:
the Mellin-Barnes representation.
JETP Lett {\bf 71}, 235-238 (2000)

\bibitem{KL01}
Kharchev, S., Lebedev, D.:
Integral representations for the eigenfunctions
of quantum open and periodic Toda chains
from the QISM formalism.
J. Phys. A: Math. Gen.{\bf 34}, 2247-2258 (2001)

\bibitem{Kos77}
Kostant, B.:
Quantisation and representation theory.
In: Representation Theory of Lie Groups,
Proc. SRC/LMS Research Symposium, Oxford 1977,
LMS Lecture Notes 34, pp. 287-316, 
Cambridge University Press, Cambridge (1977)

\bibitem{Leb65}
Lebedev, N. N.:
Special Functions and Their Applications.
Prentice-Hall, Inc. (1965)

\bibitem{MY00}
Matsumoto, H., Yor, M.:
An analogue of Pitman's $2M-X$ theorem
for exponential Wiener functionals,
Part I: A time-inversion approach.
Nagoya Math. J. {\bf 159}, 125-166 (2000)

\bibitem{MY05}
Matsumoto, H., Yor, M.:
Exponential functionals of Brownian motion.
I: Probability laws at fixed time.
Probab. Surveys {\bf 2}, 312-347 (2005)

\bibitem{OCo09}
O'Connell, N.:
Directed polymers and the quantum Toda lattice.
to appear in Ann. Probab. 
{\sf arXiv:math.PR/0910.0069}

\bibitem{OCo12}
O'Connell, N.:
Whittaker functions and related stochastic processes.\\
{\sf arXiv:math.PR/1201.4849}

\bibitem{Pit75}
Pitman, J. W.:
One-dimensional Brownian motion and the three-dimensional
Bessel process.
Adv. Appl. Prob. {\bf 7}, 511-526 (1975)

\bibitem{STS94}
Semenov-Tian-Shansky, M.:
Quantisation of open Toda lattices.
In: Dynamical Systems VII : 
Integrable Systems, Nonholonomic Dynamical Systems.
Edited by V. I. Arnol'd and S. P. Novikov.
Encyclopedia of Mathematical Sciences,
vol.16, Springer, Berlin (1994).

\bibitem{Skl85}
Sklyanin, E. K.:
The quantum Toda chain.
In: Non-linear Equations in Classical
and Quantum Field Theory, 
Lect. Notes in Physics,
{\bf 226}, pp. 195-233, Springer, Berlin (1985)

\bibitem{Sta90}
Stade, E.:
On explicit integral formulas for
$GL(n,\R)$-Whittaker functions.
Duke Math. J. {\bf 60}, 313-362 (1990)

\bibitem{Sta01}
Stade, E.:
Mellin transforms of $GL(n,\R)$ Whittaker functions.
Amer. J. Math. {\bf 123}, 121-161 (2001)

\bibitem{Stanley99}
Stanley, R. P.:
Enumerative Combinatorics. vol.2, 
Cambridge University Press, Cambridge (1999)

\bibitem{Toda89}
Toda, M.:
Theory of Nonlinear Lattices. 2nd edn. 
Springer, Berlin (1989) 

\bibitem{Wat44}
Watson, G. N.:
A Treatise on the Theory of Bessel Functions. 2nd edn.
Cambridge University Press, Cambridge (1944)

\end{thebibliography}
\end{document}